\documentclass[12pt,reqno]{amsart}
\usepackage[margin=1in]{geometry}
\usepackage{amsmath,amsthm,amssymb,comment,thmtools}
\usepackage{mathtools}
\usepackage{graphicx}
\usepackage{enumitem}
\usepackage{xcolor}
\usepackage{hyperref}
\hypersetup{
	colorlinks,
	citecolor=black,
	filecolor=black,
	linkcolor=black,
	urlcolor=black
}
\usepackage[capitalize, nameinlink, noabbrev, nosort]{cleveref}
\usepackage{float}
\usepackage{tikz}
\usetikzlibrary{math,calc}

\definecolor{cbo}{RGB}{230,159,0}
\definecolor{cbb}{RGB}{86,180,233}

\newtheorem{thm}{Theorem}
\newtheorem{cor}{Corollary}
\newtheorem{prop}{Proposition}
\newtheorem{lem}{Lemma}
\newtheorem{conj}{Conjecture}
\theoremstyle{definition}
\newtheorem{definition}{Definition}
\newtheorem{example}{Example}
\newtheorem{remark}{Remark}

\newcommand{\C}{\mathbb C}
\newcommand{\N}{\mathbb N}
\newcommand{\Z}{\mathbb Z}
\newcommand{\x}{\mathbf x}
\newcommand{\y}{\mathbf y}
\newcommand{\aarea}{\mathbf{area}}
\newcommand{\ddinv}{\mathbf{dinv}}
\newcommand{\maj}{\mathbf{maj}}
\newcommand{\bbounce}{\mathbf{bounce}}

\DeclareMathOperator{\dg}{dg}
\DeclareMathOperator{\PF}{PF}
\DeclareMathOperator{\area}{area}
\DeclareMathOperator{\dinv}{dinv}
\DeclareMathOperator{\Area}{Area}
\DeclareMathOperator{\Dinv}{Dinv}

\DeclareMathOperator{\bounce}{bounce}
\DeclareMathOperator{\rev}{rev}
\DeclareMathOperator{\sgn}{sgn}
\DeclareMathOperator{\rank}{rank}

\begin{document}

\title{A decomposition of $m$-Dyck paths}
\author{Yuhan Jiang}
\address{Department of Mathematics,
University of California, Berkeley, CA, USA}
\email{jyh@math.berkeley.edu}

\begin{abstract}
We construct a decomposition of an $m$-Dyck path into an $m$-tuple of Dyck paths such that the area sequence and bounce sequence of the $m$-Dyck path are entrywise the sum of the area sequences and bounce sequences of the Dyck paths in the tuple.
\end{abstract}

\maketitle

\section{Introduction}
The symmetric group $S_n$ acts diagonally on 
$\C[\mathbf{x},\mathbf{y}] := \C[x_1,y_1,\dots,x_n,y_n]$.
The \emph{diagonal coinvariant ring} is defined to be 
\[
  DR_n = \C[\x,\y]/\mathfrak{m}_+^{S_n}(\x,\y), \quad \mathfrak{m}_+^{S_n}(\x,\y) = \left\langle \sum_k x_k^i y_k^j: (i,j) \neq (0,0) \right \rangle
\]
the quotient space by the ideal generated by nonconstant diagonally symmetric polynomials. 
This space was studied by Haiman \cite{Haiman} for $n!$-conjecture \cite{hilbertscheme}, and its Frobenius character was studied as the shuffle conjecture \cite{10.1215/S0012-7094-04-12621-1,MR3787405}.

Let 
\[
  \mathcal{A} = \langle f: (\sigma f) = \sgn(\sigma) f, ~~\sigma \in S_n \rangle
\]
be the ideal generated by all \emph{alternating polynomials}. 
The $q,t$-Catalan numbers were first defined by Haiman as the \emph{bigraded Hilbert series} of $\mathcal{A}$, that is, $\sum_{i,j \geq 0} \dim(\mathcal{A}_{i,j}) q^i t^j$ where $\mathcal{A}_{i,j}$ is the bihomogeneous component of $\mathcal{A}$ in bidegree $(i,j)$, for which Haglund \cite{MR1972636} later gave combinatorial formulas in terms of Dyck path statistics.

Garsia and Haiman \cite{garsiahaiman} defined the \emph{space of generalized diagonal coinvariants} as
\[
  DR_n^{(m)} := (\mathcal{A}^{m-1}/\mathcal{A}^{m-1}\mathfrak{m}_+^{S_n}(\x,\y)) \otimes \epsilon^{\otimes (m-1)},
\]
where $\epsilon$ is the 1-dimensional sign representation.
The alternating component of $DR_n^{(m)}$ is
\[
  M^{(m)} \cong \left( \mathcal{A}^m/\langle \mathbf{x},\mathbf{y} \rangle \mathcal{A}^m \right) \otimes \epsilon^{\otimes (m-1)}.
\]
The Frobenius character of this ring was studied as the rational shuffle conjecture \cite{MR2971696,MR3239135,MR3383172,MR3556418,10.1215/S0012-7094-04-12621-1}.
The \emph{$q,t$-Fu\ss-Catalan numbers} were defined by Haiman \cite{Haiman} as the bigraded Hilbert series of the alternating component of the space of generalized diagonal coinvariants.

In \cite{MR2657714}, Stump asked if one could obtain a minimal generating set of $\mathcal{A}^m$ as a $\C[\x,\y]$-module, as their image in the quotient form a vector space basis for $M^{(m)}$ by the graded Nakayama's lemma. We construct a decomposition of rational Dyck paths and conjecture that such a basis is given by products of bivariate Vandermonde determinants coming from our decomposition in \cref{thm:main2}, \cref{sec:fuss}.

In \cref{sec:prelim} we collect the Dyck path statistics needed to state the conjecture, together with two combinatorial lemmas on the dinv sequence: the pairs $(a_i, d_i)$ of a Dyck path are pairwise distinct (\cref{lem:distinct}), so the bivariate Vandermonde determinants entering \cref{thm:main2} are nonzero, and a Dyck path is uniquely determined by its set of pairs (\cref{lem:injective}); consequently the determinants $\{\Delta_\pi\}$ are linearly independent already in $\C[\x,\y]$ (\cref{cor:independent}).

\section*{Acknowledgements}
The author would like to thank Christian Stump for suggesting the problem. We thank Christian Stump, Alexei Oblomkov, Mike Zabrocki, Fran\c{c}ois Bergeron, Sylvie Corteel, Lauren K. Williams, Nupur Jain, and Mitsuki Hanada for helpful discussions. We thank Alexei Oblomkov for his code. We thank Sylvie Corteel and Mitsuki Hanada for their comments on the writing.

\section{Preliminaries on Dyck paths}\label{sec:prelim}

A \emph{lattice path} is a sequence of north $(0,1)$ and east $(1,0)$ steps in the first quadrant of the plane, starting at the origin. We write $L_{m,n}$ for the set of lattice paths ending at $(m,n)$, and $L_{m,n}^+$ for the subset of paths staying weakly above the diagonal $y = \frac{n}{m}x$. A \emph{Dyck path} of semilength $n$ is an element of $L_{n,n}^+$.

Let $\pi \in L_{n,n}^+$. A cell $(i,j)$ of the $n \times n$ square is an \emph{area cell} of $\pi$ if it lies below $\pi$ and strictly above the diagonal. The \emph{area sequence} of $\pi$ is the tuple
\[
\aarea(\pi) = (a_1(\pi),\dots,a_n(\pi))
\]
such that $a_i(\pi)$ is the number of area cells of $\pi$ in the $i$-th row from the bottom \cite{MR2371044}. Area sequences are characterized by $a_1 = 0$ and $0 \leq a_{i+1} \leq a_i + 1$ for $1 \leq i < n$, and the area sequence determines the Dyck path.

We say that $(i,j)$ with $i < j$ is a \emph{primary dinv pair} of $\pi$ if $a_i = a_j$, and a \emph{secondary dinv pair} if $a_i = a_j + 1$. Let $\Dinv(\pi)$ denote the set of dinv pairs of $\pi$.

\begin{definition}\label{def:dinv}
The \emph{dinv sequence} of $\pi$ is the tuple $\ddinv(\pi) = (d_1(\pi),\dots,d_n(\pi))$ such that
\begin{equation}\label{eq:count}
d_i(\pi) = |\{j > i : a_j \in \{a_i,\, a_i - 1\}\}|
\end{equation}
is the number of dinv pairs of $\pi$ with first coordinate $i$. Adding up the dinv sequence recovers the dinv statistic of \cite{MR2371044}.
\end{definition}

For $X = ((\alpha_1,\beta_1),\dots,(\alpha_n,\beta_n)) \subseteq \N \times \N$, we define the \emph{bivariate Vandermonde determinant}
\[
\Delta_X(\x,\y) := \det\bigl(x_i^{\alpha_j} y_i^{\beta_j}\bigr)_{i,j=1}^n.
\]
For a Dyck path $\pi$, we set $X(\pi) := \bigl((d_i(\pi), a_i(\pi))\bigr)_{i=1}^n$ and $\Delta_\pi := \Delta_{X(\pi)}$.

\begin{lem}\label{lem:distinct}
For every Dyck path $\pi \in L_{n,n}^+$, the pairs $(a_i, d_i)$, $1 \leq i \leq n$, are pairwise distinct. Consequently $\Delta_\pi \neq 0$.
\end{lem}

\begin{proof}
Suppose $a_i = a_j$ with $i < j$. Every $k > j$ forming a dinv pair with $j$ satisfies $a_k \in \{a_j, a_j - 1\} = \{a_i, a_i - 1\}$ and hence forms a dinv pair with $i$; moreover, $(i,j)$ is itself a primary dinv pair, counted in $d_i$ but not in $d_j$. Hence $d_i \geq d_j + 1$, so the pairs are distinct. For the consequence, expand
\[
\Delta_\pi = \sum_{w \in S_n} \sgn(w)\, \prod_{i=1}^n x_{w(i)}^{d_i}\, y_{w(i)}^{a_i}.
\]
Since the pairs $(d_i, a_i)$ are distinct, the $n!$ monomials appearing are pairwise distinct, so no cancellation occurs.
\end{proof}

\begin{lem}\label{lem:injective}
Distinct Dyck paths $\pi \neq \pi'$ of semilength $n$ have distinct sets $\{(d_i, a_i)\}_{i=1}^n$. (By \cref{lem:distinct} the pairs within one path are distinct, so sets and multisets coincide.)
\end{lem}

\begin{proof}
Suppose $\pi \neq \pi'$ have equal pair sets, and let $i$ be the first position at which the pair \emph{sequences} differ; such $i$ exists since the area sequence determines the path. Since the prefixes agree and the global sets agree, the pairs at positions $\geq i$ form the same set $M^\ast$ in both paths; in particular, the multisets of areas at positions $\geq i$ agree. If $a_i(\pi) = a_i(\pi')$, then \eqref{eq:count} computes both $d_i$'s from that common area multiset, forcing equal pairs at position $i$, a contradiction. Hence $c := a_i(\pi) \neq c' := a_i(\pi')$.

Now $\pi$ carries the pair $P = (c, d)$ at position $i$, and $\pi'$ carries $P$ as well, at some position $j > i$. Applying \eqref{eq:count} to $P$ at both locations, the counts of areas $c$ and $c-1$ among positions $\geq i$ of $\pi'$ and among positions $\geq j$ of $\pi'$ must coincide; that is, \emph{no position of $\pi'$ in $[i, j-1]$ carries area $c$ or $c-1$}. In particular, $c' \notin \{c, c-1\}$. Symmetrically, writing $P' = (c', d')$ for the pair at position $i$ of $\pi'$, which $\pi$ carries at some position $j' > i$: no position of $\pi$ in $[i, j'-1]$ carries area $c'$ or $c'-1$, and $c \notin \{c', c'-1\}$. Together, $|c - c'| \geq 2$; without loss of generality, $c' \geq c + 2$.

Finally, consider $\pi$ on the positions $i, \dots, j'$. The positions $i, \dots, j'-1$ all avoid the band $\{c'-1,\, c'\}$, and $a_i = c \leq c' - 2$. Since consecutive areas rise by at most one, a walk starting at a value $\leq c'-2$ that never lands in the band stays $\leq c'-2$: from a value $v \leq c'-2$ the next value is $\leq v + 1 \leq c'-1$, and avoiding the band forces it back to $\leq c'-2$. Hence $a_{j'-1}(\pi) \leq c'-2$, contradicting $a_{j'}(\pi) = c' \leq a_{j'-1}(\pi) + 1$. (This includes the case $j' = i+1$, where $a_{j'-1} = a_i = c$.)
\end{proof}

\begin{cor}\label{cor:independent}
The bivariate Vandermonde determinants $\{\Delta_\pi\}$, over all Dyck paths $\pi$ of semilength $n$, are linearly independent in $\C[\x,\y]$.
\end{cor}

\begin{proof}
By the proof of \cref{lem:distinct}, the monomial support of $\Delta_\pi$ is the $S_n$-orbit of $\prod_i x_i^{d_i} y_i^{a_i}$, which is determined by the set $\{(d_i, a_i)\}$; by \cref{lem:injective} these supports are pairwise disjoint.
\end{proof}

A \emph{parking function} $P = (\pi, \sigma)$ of length $n$ consists of a Dyck path $\pi \in L_{n,n}^+$ together with a labelling of its rows by a permutation $\sigma \in S_n$ that is increasing along each vertical wall; that is, $\sigma_i < \sigma_{i+1}$ whenever $a_{i+1} = a_i + 1$. We call $\sigma$ the \emph{labelling permutation} and denote the set of parking functions of length $n$ by $\PF(n)$. The integer vector $\aarea(P)_{\sigma^{-1}}$, whose $\sigma_i$-th entry is $a_i$, is called the \emph{major sequence} of $P$, denoted by $\maj(P)$.

\begin{definition}[lex-labelling]\label{def:incrow}
We define a map $\phi: L_{n,n}^+ \to \PF(n)$ as follows. Given a Dyck path $\pi$, we order its rows in increasing order of the area; rows with the same area are ordered in increasing order of the row index. Then we label the vertical steps of $\pi$ from $1$ through $n$ with respect to this ordering.
\end{definition}

\section{$q,t$-Fu\ss-Catalan combinatorics}\label{sec:fuss}
We conjecture a basis for $\mathcal{A}^m/\langle\x,\y\rangle\mathcal{A}^m$ labeled by \emph{$m$-Dyck paths}.
An $m$-Dyck path is a lattice path from $(0,0)$ to $(mn,n)$ that stays weakly above the diagonal $my = x$.
In particular, we introduce the \emph{bounce sequence} and define a decomposition of $m$-Dyck paths into an $m$-tuple of Dyck paths which is additive on both the area sequence and the bounce sequence.

\subsection{Bounce}
To state our conjectural basis, we switch from (area,dinv) to (bounce,area).

We first introduce the \emph{bounce} defined by Haglund on Dyck paths. To define the bounce statistic, we need to define the \emph{bounce path}.

\begin{definition}[\cite{MR2371044}]
    Given a Dyck path $\pi \in L_{n,n}^+$, define the \emph{bounce path} of $\pi$ to be the path described by the following algorithm.

    Start at $(0,0)$ and travel north along $\pi$ until you encounter the beginning of an east step. Then turn east and travel until you hit the diagonal. Then turn north again, etc. Continue this way until you arrive at $(n,n)$.
    The bounce path will hit the diagonal at places $(0,0),(j_0,j_0),(j_1,j_1),\dots,(j_{b-1},j_{b-1}),(n,n)$. The $j_i$'s are called \emph{touch points}.
    We define the \emph{bounce statistic} to be the sum
    \[
    \bounce(\pi) = \sum_{i=0}^{b-1} n-j_i
    \]
\end{definition}

The idea of \emph{bounce sequence} is implicit in \cite{MR3787405}.
\begin{definition}
Given a Dyck path $\pi \in L_{n,n}^+$ whose bounce path has touch points $j_0,\dots,j_{b}$, the \emph{bounce sequence} of $\pi$ is an increasing $n$-tuple $\bbounce(\pi)$ with $j_0$ many zeroes, $j_1-j_0$ many ones, \dots, and $(n-j_{b-1})$ many $b$'s.
The sequence $(j_0, j_1 - j_0, \dots, n-j_{b-1})$ of vertical steps in the bounce path are also called the \emph{bounce composition} of $\pi$.
\end{definition}

\begin{figure}
\begin{tikzpicture}[scale=0.5]
\draw[dotted] (0,0) grid (4,4);
\draw[dotted] (0,0)--(4,4);
\draw[thick] (0,0)--(0,1)--(0,2)--(1,2)--(1,3)--(2,3)--(3,3)--(3,4)--(4,4);
\end{tikzpicture}
\hspace*{1cm}
\begin{tikzpicture}[scale=.5]
\draw[dotted] (0,0) grid (4,4);
\draw[dotted] (0,0)--(4,4);
\draw[thick] (0,0)--(0,2)--(2,2)--(2,3)--(3,3)--(3,4)--(4,4);
\node at (.5,.5) {0};
\node at (1.5,1.5) {0};
\node at (2.5,2.5) {1};
\node at (3.5,3.5) {2};
\end{tikzpicture}
\caption{A Dyck path on the left with bounce path on the right, with bounce sequence $(0,0,1,2)$.}
\label{fig:bounce0012}
\end{figure}

\begin{example}
In \cref{fig:bounce0012}, we draw a Dyck path and its bounce path, with bounce sequence labelling the diagonal blocks. This labelling scheme appeared in \cite{MR3787405}.
\end{example}

Loehr extended the bounce statistic to $m$-Dyck paths.
\begin{definition}[\cite{MR2134172}]
    Given an $m$-Dyck path $\pi \in L_{mn,n}^+$, define the \emph{bounce path} of $\pi$ to be an alternating sequence of vertical moves of lengths $v_0,v_1,\dots,v_s$ and horizontal moves of lengths $h_0,h_1,\dots$. These lengths are calculated as follows.

    Start at $(0,0)$ and travel north along $\pi$ until you encounter the beginning of an east step; the distance traveled is $v_0$.
    Then travel $v_0$ units east, so $h_0 = v_0$.
    Next, travel $v_1$ units north until you encounter an east step, and travel $v_0+v_1$ units east.
    In general, we always travel north $v_i$ units north until we are blocked by the path.
    Afterwards, for $i < m$, we travel $h_i = v_0+v_1+\cdots+v_i$ units east;
    for $i \geq m$, we travel $h_i = v_i + v_{i-1} + \cdots + v_{i-(m-1)}$ units east.
    Finally, we define the \emph{bounce statistic} to be
    \begin{equation}
      \bounce(\pi) = \sum_{k\geq 0} k v_k = \sum_{k=0}^{s} (n-v_0-v_1-\cdots-v_k),
    \end{equation}
    a weighted sum of the lengths of the vertical segments in the bounce path.
    The vertical steps $(v_0,\dots,v_s)$ are also called the \emph{bounce composition} of $\pi$.

    Conversely, a vector $(v_0,\dots,v_s)$ is a bounce composition for an $m$-Dyck path only if
    \begin{enumerate}
      \item $v_0 > 0$, $v_s > 0$;
      \item $\sum v_i = n$;
      \item there is never a string of $m$ or more consecutive zeroes.
    \end{enumerate}
\end{definition}

We generalize our definition of \emph{bounce sequence} to $m$-Dyck paths.
\begin{definition}
Given an $m$-Dyck path $\pi \in L_{mn,n}^+$, with bounce composition $(v_0,\dots, v_s)$, we define its bounce sequence to be the increasing $n$-tuple $\bbounce(\pi)$ with $v_0$ many zeroes, $v_1$ many ones, \dots, and $v_s$ many $s$'s.
The bounce sequence sums up to the bounce statistic.
\end{definition}

\begin{figure}[H]
\begin{tikzpicture}[scale=.5]
\draw[dotted] (0,0) grid (6,3);
\draw[thick] (0,0)--(0,1)--(1,1)--(1,2)--(2,2)--(3,2)--(3,3)--(4,3)--(5,3)--(6,3);
\end{tikzpicture}
\caption{A rational Dyck path in $L_{mn,m}^+$ for $m=2, n=3$, with bounce sequence (0,1,2).}
\label{fig:111}
\end{figure}

\begin{example}\label{ex:111}
The Dyck path in \cref{fig:111} is equal to its own bounce path, with $v_0 = v_1 = v_2 = 1$ and bounce sequence $(0,1,2)$. We have $h_0 = v_0 = 1, h_1 = v_0 + v_1 = 2, h_2 = v_1 + v_2 = 2, h_3 = v_2 = 1$.
\end{example}

\subsection{The decomposition algorithm}

\begin{lem}\label{lem:hv}
For any rational Dyck path $\pi \in L_{mn,n}^+$ for which the bounce path has vertical steps $v_0,v_1,\dots,v_{s}$ and horizontal steps $h_0,h_1,\dots$ (the tail terms are conventionally set to 0), for any $i \in \{0,1,\dots,m-1\}$, and for any $r \geq 0$, we have that
\begin{equation}
\sum_{k = 0}^{r} h_{km+i} = \sum_{j=0}^{rm+i} v_j.
\end{equation}
In particular, $\sum_{k \geq 0} h_{km+i} = n$.
\end{lem}
\begin{proof}
By definition, for $j < m$, we have 
\[h_j = v_j+\cdots+v_1+v_0.\]
For $j \geq m$, we have
\[
  h_j = v_j + v_{j-1} + \cdots + v_{j-m+1}.
\]
Therefore, 
\begin{equation}
  \sum_{k=0}^r h_{km+i} = (v_i + \cdots + v_0) + \sum_{k=0}^{r-1} (v_{km+m+i}+\cdots+v_{km+i+1}) = v_0 + \cdots + v_{rm+i}.
\end{equation}
Setting $r$ to be sufficiently large, we have that $\sum_{k\geq 0} h_{km+i} = \sum v_i = n$. 
\end{proof}

\begin{definition}\label{def:decom}
Let $\pi \in L_{mn,n}^+$ be an $m$-Dyck path of which the bounce path has horizontal steps $h_0,h_1,\dots$.
For any $i \in \{0,1,\dots,m-1\}$, we define $\pi^i \in L_{n,n}^+$ to be the concatenation of the parts under $\pi$ in the columns specified by $h_{km+i}$ for $k\geq 0$.
In particular, the bounce composition of $\pi^i$ is given by $(h_i,h_{m+i},\dots)$.
\end{definition}

\begin{figure}[H]
\begin{tikzpicture}[scale=.5]
\draw[dotted] (0,0) grid (6,3);
\draw[thick] (0,0)--(0,1)--(1,1)--(1,2)--(2,2)--(2,3)--(3,3)--(4,3)--(5,3)--(6,3);
\draw[fill=cbb] (1,0)--(3,0)--(3,3)--(2,3)--(2,2)--(1,2)--(1,0);
\draw[fill=cbb] (5,0)--(6,0)--(6,3)--(5,3)--(5,0);
\draw[fill=cbo] (0,0)--(1,0)--(1,1)--(0,1)--(0,0);
\draw[fill=cbo] (3,0)--(5,0)--(5,3)--(3,3)--(3,0);
\draw node at (2,-.5) {area: (0,1,2)};
\draw node at (2.5,-1.3) {bounce: (0,1,2)};
\end{tikzpicture}
\hspace{1cm}
\begin{tikzpicture}[scale=.5]
\draw[dotted] (0,0) grid (3,3);
\draw[thick] (0,0)--(0,1)--(1,1)--(1,2)--(1,3)--(2,3)--(3,3);
\draw node at (1,-.5) {area: (0,0,1)};
\draw node at (1.5,-1.3) {bounce: (0,1,1)};
\end{tikzpicture}
\begin{tikzpicture}[scale=.5]
\draw[dotted] (0,0) grid (3,3);
\draw[thick] (0,0)--(0,1)--(0,2)--(1,2)--(1,3)--(2,3)--(3,3);
\draw node at (1,-.5) {area: (0,1,1)};
\draw node at (1.5,-1.3) {bounce: (0,0,1)};
\end{tikzpicture}
\caption{A rational Dyck path (left) with bounce path shown in \cref{fig:111}, bounce sequence $(0,1,2)$, and area sequence $(0,1,2)$, and its decomposition into two Dyck paths (right).}
\label{fig:decomp}
\end{figure}

\begin{example}
The 2-Dyck path in \cref{fig:111} is the bounce path of the 2-Dyck path on the left in \cref{fig:decomp}. In particular, we have the same horizontal steps $(1,2,2,1)$ as in \cref{ex:111}. We decompose the 2-Dyck path with respect to the horizontal steps and color them alternatingly. Then we collect regions of the same color and concatenate them into two Dyck paths on the right.
\end{example}

\begin{thm}\label{thm:decomp}
For any $m$-Dyck path $\pi \in L_{mn,n}^+$ and for any $i \in \{0,1,\dots,m-1\}$, let $\pi^i$ be defined as in \cref{def:decom}.
Then, the bounce (resp. area) sequence of $\pi$ is equal to the entrywise sum of the bounce (resp. area) sequences of $\pi^i$, i.e.,
\begin{align}
\aarea(\pi) &= \sum_{i=0}^{m-1} \aarea(\pi^i) \\
\bbounce(\pi) &= \sum_{i=0}^{m-1} \bbounce(\pi^i).
\end{align}
\end{thm}

To prove this theorem, we need the following definitions

\begin{definition}
A \emph{partition} $\lambda = (\lambda_1\geq \lambda_2 \geq \cdots \geq \lambda_n \geq 0)$ is a weakly decreasing sequence of nonnegative integers.
The \emph{conjugate partition} of $\lambda$ is denoted as $\lambda'$ where $\lambda'_j = |\{i: \lambda_i \geq j\}|$.
\end{definition}

\begin{definition}
Let $\pi \in L_{n,n}^+$ be a Dyck path with area sequence $\aarea(\pi) = (a_1,\dots,a_n)$.
The \emph{diagram} of $\pi$ is the partition $\lambda$ such that $\lambda_i = n-i-a_{n-i+1}$.

Let $\pi \in L_{mn,n}^+$ be an $m$-Dyck path with area sequence $\aarea(\pi) = (a_1,\dots,a_n)$.
The \emph{diagram} of $\pi$ is the partition $\lambda$ such that $\lambda_i = m(n-i)-a_{n-i+1}$.
\end{definition}

\begin{proof}[Proof of \cref{thm:decomp}]
By \cref{def:decom}, we see that for any $k \in [mn]$, there exists $i \in \{0,\dots,m-1\}$ and $j \in [n]$ such that $\dg(\pi)'_k = \dg(\pi^i)'_j$.
Conversely, for any $i \in [0,m-1]$ and any $j \in [n]$, there exists $k \in [mn]$ such that $\dg(\pi^i)'_j = \dg(\pi)'_k$.
In other words, $\dg(\pi)$ is a concatenation of the parts of $\pi^i$s.
Then, $\dg(\pi)$ as the conjugate of $\dg(\pi)'$, satisfies
$\dg(\pi)_\ell = |\{k: \dg(\pi)'_k \geq \ell\}| = \sum |\{k: \dg(\pi^i)'_j \geq \ell\}| = \sum_{i=0}^{m-1} \dg(\pi^i)_\ell$.
As the area sequence is determined by the diagram, we have the first equality.

For the second equality, observe that for a weakly increasing $n$-tuple with $c_r$ many entries equal to $r$ for each $r \geq 0$, the $k$-th entry equals $\#\{r \geq 0 : c_0 + c_1 + \cdots + c_r < k\}$.
By \cref{def:decom}, the bounce composition of $\pi^i$ is $(h_i, h_{m+i}, h_{2m+i}, \dots)$, so $\bbounce(\pi^i)$ has $h_{rm+i}$ many entries equal to $r$. Applying the observation and then \cref{lem:hv},
\[
\bbounce(\pi^i)_k
\;=\; \#\Bigl\{r \geq 0 : \sum_{\ell=0}^{r} h_{\ell m + i} < k\Bigr\}
\;=\; \#\Bigl\{r \geq 0 : \sum_{j=0}^{rm+i} v_j < k\Bigr\}.
\]
Summing over $i \in \{0,\dots,m-1\}$ and using that $(i,r) \mapsto j = rm+i$ is a bijection from $\{0,\dots,m-1\} \times \Z_{\geq 0}$ onto $\Z_{\geq 0}$,
\[
\sum_{i=0}^{m-1} \bbounce(\pi^i)_k
\;=\; \#\Bigl\{j \geq 0 : \sum_{\ell=0}^{j} v_\ell < k\Bigr\}
\;=\; \bbounce(\pi)_k,
\]
where the last equality is the same observation applied to $\bbounce(\pi)$, which has $v_j$ many entries equal to $j$.
\end{proof}

\begin{remark}
Athanasiadis \cite{MR2098091} showed a way to decompose an $m$-Dyck path into an $m$-tuple of Dyck paths for which the area sequence is additive. Our decomposition is essentially different.
\end{remark}

\subsection{The sweep map}
Haglund introduced the $\zeta$-map for Dyck paths\cite{MR2163448}, and then Thomas and Williams generalized it to rational Dyck paths and even broader context \cite{MR2371044}. 
We follow the conventions in \cite{MR3787405}.
It is implied by classical results that the $\zeta$-map, combined with the \cref{def:incrow} lex-labelling of Dyck paths, is a bijection between the bounce sequences, area sequences, and dinv sequences.

An $(m,n)$-Dyck path is a lattice path from $(0,0)$ to $(m,n)$ staying weakly above the diagonal $y = \frac{n}{m}x$. 
Thomas--Williams definition applies to $(m,n)$-Dyck paths for arbitrary positive integers $m,n$.
\begin{definition}\label{def:rank}
Given an $(m,n)$-Dyck path $\pi$, 
let $m' = m/\gcd(m,n)$ and $n' = n/\gcd(m,n)$.
We define the \emph{rank} of a lattice point $(x,y)$ in the $m \times n$-rectangle as
$\rank(x,y) = m' y - n' x.$
Let the rank of a vertical step in $\pi$ be the rank of its south end, and let the rank of a horizontal step be the rank of its west end.
Then $\zeta(\pi)$ is obtained as we reorder the steps in $\pi$ by increasing rank, and order the steps of the same rank from northeast to southwest.
\end{definition}

\begin{figure}[H]
\centering
\begin{tikzpicture}[scale=0.5]
\draw[dotted] (0,0) grid (6,4);
\draw[thick] (0,0)--(0,1)--(0,2)--(0,3)--(1,3)--(2,3)--(3,3)--(3,4)--(4,4)--(5,4)--(6,4);
\end{tikzpicture} 
\begin{tikzpicture}[scale=0.5]
\draw[dotted] (0,0) grid (6,4);
\draw[thick] (0,0)--(0,1)--(1,1)--(1,2)--(1,3)--(2,3)--(3,3)--(4,3)--(4,4)--(5,4)--(6,4);
\end{tikzpicture}
\caption{A (6,4)-Dyck path on the left and its $\zeta$-image on the right.}
\label{fig:64}
\end{figure}

\begin{example}
A Dyck word of a Dyck path $\pi \in L_{m,n}^+$ is a binary word in $\{0,1\}^{m+n}$ where 1 represents a north step and 0 represents an east step.
Consider the (6,4)-Dyck path with Dyck word $w = (1,1,1,0,0,0,1,0,0,0)$ and rank displayed in the second row of the array below.
We obtain $w' = (1,0,1,1,0,0,0,1,0,0)$ as we order the elements based on their ranks.
\[
  \begin{pmatrix}
  1 & 1 & 1 & 0 & 0 & 0 & 1 & 0 & 0 & 0 \\
  0 & 3 & 6 & 9 & 7 & 5 & 3 & 6 & 4 & 2
  \end{pmatrix}
\]
\end{example}

\begin{thm}[\cite{MR1972636,MR2163448,ThomasWilliams2018}]
For any Dyck path $\pi$, we have $\dinv(\pi) = \area(\zeta(\pi))$ and $\area(\pi) = \bounce(\zeta(\pi))$.
\end{thm}

\begin{remark}
The reversal of a Dyck path $\pi$ is obtained by flipping it across the diagonal, denoted $\rev(\pi)$. The convention we follow here results in the reversal of the $\zeta$-map defined in \cite{MR2371044} and in \texttt{SageMath}.
\end{remark}

Carlsson and Mellit gave an alternative description for the $\zeta$-map of Dyck paths.
Let $\Dinv(\pi)$ be the set of dinv pairs of $\pi$ and let $\Area(\pi)$ be the set of area cells of $\pi$.
Let $\sigma$ be the labelling permutation of the parking function $\phi(\pi)$.
Then $\zeta(\pi)$ is defined by 
\begin{equation}
  \Area(\zeta(\pi)) = \{(\sigma_j,\sigma_{j'}): (j,j') \in \Dinv(\pi)\}.
\end{equation}

Their results further implied the following lemma.
\begin{lem}\label{prop:vector_zeta}
Let $\pi \in L_{n,n}^+$ be a Dyck path such that $\phi(\pi) = (\pi,\sigma)$, that is, $\sigma$ is the labelling permutation of the parking function obtained from $\pi$ through its lex-labelling, see \cref{def:incrow}. Then $\bbounce(\zeta^{-1}(\pi)) = \maj(\phi(\pi))$, and $\aarea(\rev(\zeta^{-1}(\pi))) = \rev(\ddinv(\pi))$. 

Moreover, primary dinv pairs in $\pi$ map to area cells below the bounce path in $\zeta(\pi)$, and secondary dinv pairs in $\pi$ map to area cells above the bounce path.
\end{lem}

\begin{figure}[H]
\begin{tikzpicture}[scale=.5]
\draw[dotted] (0,0) grid (4,4);
\draw[thick] (0,0)--(0,1)--(0,2)--(1,2)--(1,3)--(2,3)--(3,3)--(3,4)--(4,4);
\end{tikzpicture}
\begin{tikzpicture}[scale=.5]
\draw[dotted] (0,0) grid (4,4);
\draw[thick] (0,0)--(0,1)--(0,2)--(1,2)--(1,3)--(2,3)--(3,3)--(3,4)--(4,4);
\draw node at (0.500000,0.500000) {1};
\draw node at (0.500000,1.500000) {3};
\draw node at (1.500000,2.500000) {4};
\draw node at (3.500000,3.500000) {2};
\end{tikzpicture}
\hspace{1cm}
\begin{tikzpicture}[scale=.5]
\draw[dotted] (0,0) grid (4,4);
\draw[dotted] (0,0) -- (4,4);
\draw[thick] (0,0)--(0,1)--(0,2)--(1,2)--(1,3)--(1,4)--(2,4)--(3,4)--(4,4);
\end{tikzpicture}
\begin{tikzpicture}[scale=.5]
\draw[dotted] (0,0) grid (4,4);
\draw[thick] (0,0)--(0,1)--(0,2)--(1,2)--(2,2)--(2,3)--(2,4)--(3,4)--(4,4);
\end{tikzpicture}
\caption{A Dyck path $\pi$ (left), its lex-labelling $\phi(\pi)$, and its image under the $\zeta$-map (right) and the bounce path of $\zeta(\pi)$.}
\label{fig:zeta}
\end{figure}

\begin{example}
Consider the Dyck path $\pi$ on the left in \cref{fig:zeta} with area sequence $\aarea(\pi) = (0,1,1,0)$.
The dinv pairs for $\phi(\pi)$ are $(1,2),(2,3),(2,4),(3,4)$, which are exactly the area cells for the Dyck path on the right $\zeta(\pi)$ in \cref{fig:zeta}.
The bounce sequence of $\zeta(\pi)$ is $\bbounce(\zeta(\pi)) = (0,0,1,1)$. The lex-labelling of $\pi$ gives the permutation $\sigma = (1,3,4,2)$, and $\sigma(\bbounce(\zeta(\pi))) = (0,1,1,0) = \aarea(\pi)$.
The area sequence of the reversal of $\zeta(\pi)$ is $\aarea(\rev(\zeta(\pi))) = (0,1,2,1)$, which is also the reversal of $\ddinv(\pi) = (1,2,1,0)$.
\end{example}

\subsection{A conjectural basis of alternating generalized diagonal coinvariants}

We conjecture that our decomposition provides a basis for the alternating component of the generalized diagonal coinvariants.
\begin{conj}\label{thm:main2}
For any rational Dyck path $\pi \in L_{mn,n}^+$, for any $i \in \{0,1,\dots,m-1\}$,  let $\pi^i \in L_{n,n}^+$ be the Dyck path defined in \cref{def:decom}.
Then $\{\prod_{i=0}^{m-1} \Delta_{\zeta^{-1}(\pi^i)}\}$ over all $m$-Dyck paths form a basis of the alternating component of the space of generalized diagonal coinvariants.
\end{conj}

\begin{figure}[H]
\centering
\begin{tikzpicture}[scale=0.5]
\draw[dotted] (0,0) grid (6,3);
\draw[thick] (0,0)--(0,1)--(1,1)--(1,2)--(2,2)--(3,2)--(3,3)--(4,3)--(5,3)--(6,3);
\draw[fill=cbo] (0,0)--(0,1)--(1,1)--(1,0);
\draw[fill=cbo] (3,0)--(3,3)--(5,3)--(5,0);
\draw[fill=cbb] (1,0)--(1,2)--(3,2)--(3,0);
\draw[fill=cbb] (5,0)--(5,3)--(6,3)--(6,0);
\draw node at (2,-.5) {area: (0,1,1)};
\draw node at (2.5,-1.2) {bounce: (0,1,2)};
\end{tikzpicture}
\hspace{.5cm}
\begin{tikzpicture}[scale=0.5]
\draw[dotted] (0,0) grid (3,3);
\draw[thick] (0,0)--(0,1)--(1,1)--(1,2)--(1,3)--(2,3)--(3,3);
\draw node at (2,-.5) {area: (0,0,1)};
\draw node at (2,-1.2) {bounce: (0,1,1)};
\end{tikzpicture}
\begin{tikzpicture}[scale=0.5]
\draw[dotted] (0,0) grid (3,3);
\draw[thick] (0,0)--(0,1)--(0,2)--(1,2)--(1,3)--(2,3)--(3,3);
\draw node at (2,-.5) {area: (0,1,1)};
\draw node at (2,-1.2) {dinv: (0,1,0)};
\end{tikzpicture}
\hspace{.5cm}
\begin{tikzpicture}[scale=0.5]
\draw[dotted] (0,0) grid (3,3);
\draw[thick] (0,0)--(0,1)--(0,2)--(1,2)--(2,2)--(2,3)--(3,3);
\draw node at (2,-.5) {area: (0,1,0)};
\draw node at (2,-1.2) {bounce: (0,0,1)};
\end{tikzpicture}
\begin{tikzpicture}[scale=0.5]
\draw[dotted] (0,0) grid (3,3);
\draw[thick] (0,0)--(0,1)--(1,1)--(1,2)--(1,3)--(2,3)--(3,3);
\draw node at (2,-.5) {area: (0,0,1)};
\draw node at (2,-1.2) {dinv: (1,0,0)};
\end{tikzpicture}
\caption{The decomposition of a $2$-Dyck path and the $\zeta^{-1}$-image of this decomposition.}
\label{fig:product}
\end{figure}

\begin{example}
In \cref{fig:product}, we show the decomposition of a $2$-Dyck path. The corresponding $\Delta_{\zeta^{-1}(\pi^1)} \Delta_{\zeta^{-1}(\pi^2)}$ is equal to
\[ 
(-x_{1} y_{1} y_{2} + x_{1} y_{1} y_{3} + x_{2} y_{1} y_{2} - x_{2} y_{2} y_{3} - x_{3} y_{1} y_{3} + x_{3} y_{2} y_{3})(x_{1} y_{2} - x_{1} y_{3} - x_{2} y_{1} + x_{2} y_{3} + x_{3} y_{1} - x_{3} y_{2}).
\]
\end{example}

\begin{definition}
An $m$-parking function consists of an $m$-Dyck path $\pi$ and a permutation $\sigma \in S_n$ which is increasing along each vertical wall.
\end{definition}
We analyze the synergy between our decomposition algorithm \cref{def:decom} and the $\zeta$-map on parking functions.
\begin{prop}
Let $(\pi,\sigma)$ be an $m$-parking function, and let $(r_1,\dots,r_n)$ be the permuted rank vector where $r_i$ is the rank of the vertical step in $\pi$ with label $i$ (see \cref{def:rank}).
Let $(v_0,v_1,\dots)$ be the bounce composition of $\zeta(\pi)$.
Then the labels of $\zeta(\pi)$ to the left of the $v_j$ vertical steps in the bounce path of $\zeta(\pi)$ are exactly those labels $i$ with $r_i = j$.
\end{prop}
\begin{proof}
The most northeast step of rank $j$ is always horizontal, because the path eventually returns to rank 0.
Therefore, we want to show that the number of horizontal steps of rank $j$ is equal to $h_{j-1}$.
We prove by induction.
The base case asserts that the number of horizontal steps of rank 1 is equal to $h_0 = v_0$.
A horizontal step has rank 1 if and only if its west end has rank 1 and its east end has rank 0.
Either its east end is $(mn,n)$, or it is followed by some vertical step of rank 0.
There is always a vertical step of rank 0 at $(0,0)$.
Therefore, we have proven the base case.

For the induction step, a horizontal step has rank $j$ if and only if its west end has rank $j$ and its east end has rank $j-1$, which can be followed by either vertical steps or horizontal steps of rank $j-1$.
However, for every vertical step of rank $j-m-1$, the north end of such a vertical step has rank $j-1$, and must be followed by some step of rank $j-1$.
Therefore, the number of horizontal steps of rank $j$ is equal to, by induction hypothesis, $v_{j-1}+h_{j-2}-v_{j-m-1} = h_{j-1}$.
\end{proof}

\bibliographystyle{alpha}
\bibliography{main}

\end{document}